\documentclass[a4paper,12pt]{amsart}
\newif\ifdraft
 \draftfalse

\usepackage[english]{babel}
\usepackage[utf8]{inputenc}
\allowdisplaybreaks  % permite cortes de página en fórmulas

\usepackage{amssymb}
\usepackage{xypic,amscd}

\usepackage{mathtools}
\usepackage{etoolbox}

\usepackage[colorlinks,linkcolor=Red3,citecolor=Red3,urlcolor=Cyan3]{hyperref}

\usepackage{enumitem}

\usepackage{csquotes}

\usepackage{parskip}

\usepackage{bm}

\newtheorem{definition}[equation]{Definition}
\newtheorem{theorem}[equation]{Theorem}

\newtheorem{corollary}[equation]{Corollary}
\newtheorem{lemma}[equation]{Lemma}
\newtheorem{proposition}[equation]{Proposition}

\theoremstyle{remark}

\theoremstyle{remark}

\numberwithin{equation}{section} %\numberwithin{defn}{section}

\renewcommand{\cH}{{\mathcal{H}}}

\renewcommand{\cR}{{\mathcal{R}}}
\newcommand{\cS}{{\mathcal{S}}}

\newcommand{\bC}{{\mathbb{C}}}

\newcommand{\bN}{{\mathbb{N}}}

\newcommand{\bQ}{{\mathbb{Q}}}
\newcommand{\bR}{{\mathbb{R}}}
\newcommand{\bZ}{{\mathbb{Z}}}

\newcommand{\gal}{\operatorname{Gal}}

\newcommand{\norm}{{\operatorname{N}}}

\newcommand{\rk}{{\operatorname{rk}}}
\usepackage[svgnames,x11names,dvipsnames]{xcolor}
\ifdraft
\RequirePackage{soulutf8}
\sethlcolor{MistyRose}
\usepackage{xparse}
\NewDocumentCommand{\badref}{v}{\textcolor{Firebrick1}{\framebox{\textbf{??}}\hl{\texttt{#1}}}}
\newtheoremstyle{comentario}%
{}%
{}%
{\sffamily}%
{0pt}%
{\upshape}%
{}%
{0pt}%
{\color{Red2}{\sffamily\thmname{#1}\framebox{\small\textbf{\thmnumber{#2}}}}\mdseries\upshape\thmnote{ (#3)} }%
\theoremstyle{comentario}

\theoremstyle{plain}
\usepackage[color,notcite,notref]{showkeys}%
\definecolor{labelkey}{rgb}{1.0,0.8,0.8}%
\usepackage[hyperpageref]{backref} % For drafts
\renewcommand*{\backrefalt}[4]{%
\ifcase#1%
\textcolor{Red1}{(No citations)}%
\or%
\textcolor{Yellow3}{(One citation on page #2.)}%
\else%
\textcolor{Green3}{(#1 citations in pages #2.)}%
\fi%
}%
\fi
\title[]{A connection between the poles of the zeta function of a recurrence sequence and the module of relations of its roots}

\author[Á. Serrano~Holgado]{Álvaro Serrano Holgado}
\email{Alvaro\_Serrano@usal.es}

\thanks{Research of the author is supported by grant PID2021-124332NB-C22 of the MICINN (Spain).}
\address{Departamento de Matemáticas, Universidad de
Salamanca,  Plaza de la Merced 1-4
        \\
        37008 Salamanca. Spain. 
}



\begin{document}

\maketitle

\begin{abstract}
Answering a question left open in previous research, we study the enumeration of poles of the zeta function $\varphi(s)$ associated to an integer linear recurrence sequence $\{a_n\}$. This enumeration can count poles more than once, and we prove that this happens if and only if the module of relations of the roots of the recurrence is nontrivial. A review of the existing literature on the module of relations yields a series of sufficient conditions for the enumeration of poles of $\varphi(s)$ to be injective. All of this is illustrated by examples of both cases.
\end{abstract}

%%%%%%%%%%%%%%%%%%%%%%%%%%%%%%%%%%%%%%%%%%%%%%%%%%%%%%%%%%%%%%%%%%%%%%%%%%%%%%%%%%%
\section{Introduction}
\label{sec:intro}
%%%%%%%%%%%%%%%%%%%%%%%%%%%%%%%%%%%%%%%%%%%%%%%%%%%%%%%%%%%%%%%%%%%%%%%%%%%%%%%%%%%
%%%%%%%%%%%%%%%%%%%%%%%%%%%%%%%%%%%%%%%

The study of the zeta functions associated to an integer linear recurrence sequence was started independently by Navas in \cite{Navas} and Egami in \cite{Egami}, both concerning the Fibonacci sequence. Since then, quite a few different aspects and generalizations of these studies have been made: see, for example, Kamano's \cite{Kamano}, dealing with the zeta and L-functions of a Lucas sequence, Meher and Rout's \cite{MeherRout1} and \cite{MeherRout2}, about multiple Lucas zeta functions. All these articles, however, study only recurrences of degree $2$.

Not until recently has the higher degree case been explored, by Smajlovi\'c et al. in \cite{SmajlovicTribonacci} and in our own papers \cite{Serrano-Tribonacci} and \cite{Serrano}. This higher degree case presents a difficulty that is not present in degree $2$: analytic continuation of the corresponding Dirichlet series is proved through representation by a different series (see Section~\ref{sec:background} below, and particularly Equation~\eqref{eq:phi}). Each term of this series has a pole, and while in the degree $2$ case these poles are all at different points, that need not happen in higher degrees. A question equivalent to this was left open in \cite[Remark~4.2]{SmajlovicTribonacci}:

\begin{displayquote}
Poles at other negative integers might exist only in case that the number [...] is some root of unity. We find this outcome highly unlikely, in view of properties of the conjugate Pisot numbers elaborated in [...] and [...]. However, we were not able to prove that the argument of a conjugate of the third order Pisot number is not a rational multiple of $\pi$.
\end{displayquote}

This particular question was answered in our article \cite[Section~5, example~2]{Serrano}, but still not in the general case: in \cite[Remark~4.6]{Serrano}, we wrote:

\begin{displayquote}
Note that we have not said that all the points [...] are distinct. [...] Though in general we do not know the answer to this question [...].
\end{displayquote}

This paper completes our previous work \cite{Serrano} and gives a characterization of the distinctness of the poles of the terms of the zeta function associated to a recurrence sequence in terms of the module of relations of the roots of said recurrence, in Theorem~\ref{t:equivalence}. The article is structured as follows:

In Section~\ref{sec:background}, we provide an overview of our previous work about zeta functions of recurrence sequences and define, in \eqref{eq:indexation}, the ``indexation map'' $\iota$ whose fibres we are interested in. In Section~\ref{sec:finite-fibres} it is proved that this map has finite fibres.

In Section~\ref{sec:equivalence}, we introduce the module of relations and we prove the main result of this paper, Theorem~\ref{t:equivalence} (or its equivalent statement, Theorem~\ref{t:equivalence-triviality}), which states that $\iota$ is one-to-one if and only if the module of relations is trivial.

In Section~\ref{sec:triviality-so-far}, we review the existing literature and results on the module of relations, and in Theorem~\ref{t:sufficient-conditions} use those results to give a few sufficient conditions for its triviality.

Finally, in Section~\ref{sec:nontriviality}, we give a few examples of numbers with nontrivial module, and in Proposition~\ref{p:fibre-at-0} we prove that the fibre of $\iota$ at the pole $s=0$ consists always of exactly one point, whether or not the module is trivial.

%%%%%%%%%%%%%%%%%%%%%%%%%%%%%%%%%%%%%%%%%%%%%%%%%%%%%%%%%%%%%%%%%%%%%%%%%%%%%%%%%%%
\section{Background}
\label{sec:background}
%%%%%%%%%%%%%%%%%%%%%%%%%%%%%%%%%%%%%%%%%%%%%%%%%%%%%%%%%%%%%%%%%%%%%%%%%%%%%%%%%%%
%%%%%%%%%%%%%%%%%%%%%%%%%%%%%%%%%%%%%%%

Given a Perron number $\alpha$ (that is, a real algebraic integer greater than $1$ such that all of its conjugates are strictly smaller than itself in modulus) of degree $r$ with conjugates $\alpha_1=\alpha, \alpha_2, \dots, \alpha_r$, and an integer linear recurrence sequence $\{a_n\}_{n\in\bN}$ satisfying the minimal polynomial $p_{\alpha}(T)\in\bZ[T]$ of $\alpha$, there are $r$ numbers $\lambda_i\in\bQ(\alpha_i)$, $1\leq i\leq r$ (see \cite[\S~1.1.6]{Everest} for more details) such that
\begin{equation}\label{eq:binet}
	a_n=\lambda_1\alpha_1^n+\dots+\lambda_r\alpha_r^n.
\end{equation}

If $\lambda_1>0$, $\{a_n\}$ is eventually strictly increasing and positive. Discarding, if necessary, a finite number of terms, we can define the Dirichlet series associated to $\{a_n\}$,
\[
	\sum_{n=1}^{\infty}\frac{1}{a_n^s}.
\]

From formula \eqref{eq:binet}, and using the binomial theorem, we proved in \cite{Serrano} that this series has an analytic continuation to a meromorphic function of the whole $s$-plane, given by the expression (a bit different from the one we gave in \cite{Serrano}, but easily seen to be equivalent)
\begin{equation}\label{eq:phi}
	\varphi(s)=\lambda_1^{-s}\sum_{\kappa\in\bN_0^{r-1}}\binom{-s}{\kappa}\left(\frac{\lambda_2}{\lambda_1},\dots,\frac{\lambda_r}{\lambda_1}\right)^{\kappa}\frac{1}{\alpha_1^{s+|\kappa|}(\alpha_2,\dots,\alpha_r)^{\kappa}},
\end{equation}
where we use the following notation:
\begin{enumerate}[label=$\bullet$]
\item $\kappa=(\kappa_2,\dots,\kappa_r)\in\bN_0^{r-1}$ is a multi-index of length $r-1$ (we number them from $2$ to $r$, instead of from $1$ to $r-1$, to make them correspond with the roots $\alpha_2,\dots,\alpha_r$).

\item Given a vector $(z_2,\dots, z_r)\in\bC^{r-1}$, $(z_2,\dots,z_r)^{\kappa}=z_2^{\kappa_2}\dots z_r^{\kappa_r}$.

\item $|\kappa|=\kappa_2+\dots+\kappa_r$, and $\kappa!=\kappa_2!\dots\kappa_r!$.

\item Given $z\in\bC$, the multinomial polynomial $\binom{z}{\kappa}$ is defined by
\[
	\binom{z}{\kappa}=\frac{z(z-1)\dots(z-|\kappa|+1)}{\kappa!}.
\]
This multinomial coefficient generalises both the binomial polynomial
\[
	\binom{z}{m}=\frac{z(z-1)\dots (z-m+1)}{m!},
\]
for $z\in\bC$ and $m\in\bN_0$, and the multinomial coefficient
\[
	\binom{|\kappa|}{\kappa}=\frac{|\kappa|!}{\kappa_2!\dots\kappa_r!},
\]
for $\kappa\in\bN_0^{r-1}$.
\end{enumerate}

From the series expansion \eqref{eq:phi} we see that the poles of $\varphi(s)$ lie among the points
\begin{equation}\label{eq:poles}
	s_{\kappa,n}=-|\kappa|+\frac{\log|\alpha_2^{\kappa_2}\dots\alpha_r^{\kappa_r}|}{\log\alpha}+i\frac{\arg(\alpha_2^{\kappa_2}\dots\alpha_r^{\kappa_r})}{\log\alpha}+\frac{2\pi i n}{\log\alpha},
\end{equation}
for $\kappa\in\bN_0^{r-1}$ and $n\in\bZ$.

We will explore under which conditions the points enumerated by~\eqref{eq:poles} are different.

First of all, note that for fixed $\kappa\in\bN_0^{r-1}$, $\{s_{\kappa,n}\}_{n\in\bZ}$ is a set of periodic points lying on the same vertical line, separated by intervals of length $\frac{2\pi}{\log\alpha}$. Given two $\kappa,\kappa'\in\bN_0^{r-1}$, the sets $\{s_{\kappa,n}\}$ and $\{s_{\kappa',n}\}$ are either the same or disjoint. Therefore, we only need to see what happens in a ``fundamental region'' (borrowing the name from the theory of elliptic functions, but keeping in mind the obvious distinction), which can be taken to be, for example, horizontal bands
\[
	\Omega_{\tau}=\left\{s\in\bC \:\big\vert\: \tau-\frac{\pi}{\log\alpha}<\Im(s)\leq \tau+\frac{\pi}{\log\alpha}\right\},
\]
for $\tau\in \bR$. Assuming by default that $\arg$ is the principal branch of the argument, namely, the one with range $(-\pi,\pi]$, we take the fundamental region in which the poles $s_{\kappa}=s_{\kappa,0}$ lie, that is,
\[
	\Omega_0=\left\{-\frac{\pi}{\log\alpha}<\Im(s)\leq\frac{\pi}{\log\alpha}\right\}
\]

Accordingly, we define our ``indexation map''
\begin{equation}\label{eq:indexation}
	\begin{array}{rccl}
	\iota: & \bN_0^{r-1} & \longrightarrow & \Omega_0 \\
	 & \kappa & \longmapsto & s_{\kappa}\end{array}.
\end{equation}

\begin{figure}[!ht]
\begin{center}	
\includegraphics[scale=1]{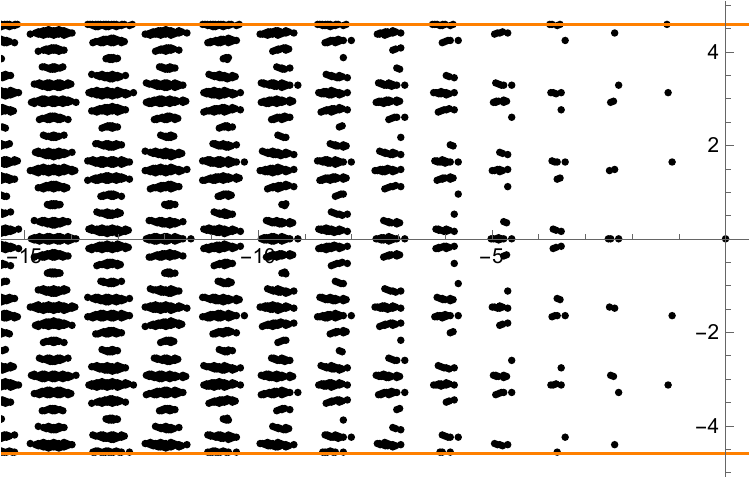}
\caption{Image of $\iota$ for $\alpha$ with minimal polynomial $P_6(x)$.}
\label{fig:hexabonacci}
\end{center}
\end{figure}

If we plot the image of $\iota$ for a couple of Perron numbers, we can see the role that the injectivity (or lack thereof) of $\iota$ plays. For the Perron number $\alpha$ which is the largest real root of $P_6(x)=x^6-x^5-x^4-x^3-x^2-x-1$, the image of $\iota$ is shown in Figure~\ref{fig:hexabonacci}.

\begin{figure}[!ht]
\begin{center}	
\includegraphics[scale=1]{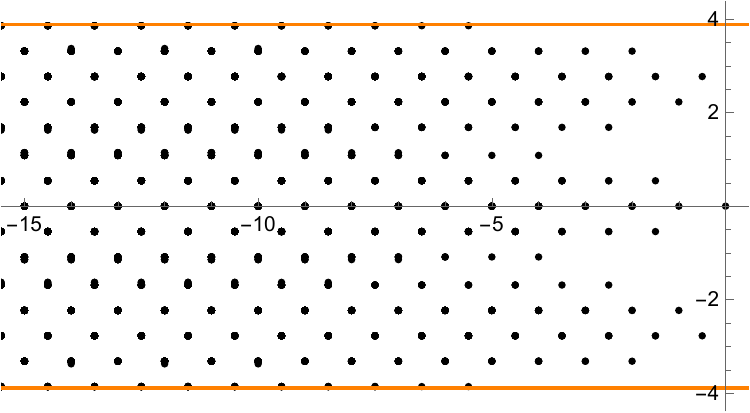}
\caption{Image of $\iota$ for $\alpha$ with minimal polynomial $f(x)$.}
\label{fig:schinzel}
\end{center}
\end{figure}

In contrast, for the Perron number $\alpha$ with minimal polynomial $f(x)=x^6-2x^4-6x^3-2x^2+1$ (which we will study further in Section~\ref{sec:nontriviality}), the image of $\iota$ is shown in Figure~\ref{fig:schinzel}.

%%%%%%%%%%%%%%%%%%%%%%%%%%%%%%%%%%%%%%%%%%%%%%%%%%%%%%%%%%%%%%%%%%%%%%%%%%%%%%%%%%%
\section{Finiteness of the fibres of $\iota$}
\label{sec:finite-fibres}
%%%%%%%%%%%%%%%%%%%%%%%%%%%%%%%%%%%%%%%%%%%%%%%%%%%%%%%%%%%%%%%%%%%%%%%%%%%%%%%%%%%
%%%%%%%%%%%%%%%%%%%%%%%%%%%%%%%%%%%%%%%

Our first result, which we stated without proof in \cite{Serrano}, is that the map $\iota$ defined in \eqref{eq:indexation} has finite fibres, that is, each pole can only be reached by at most a finite number of $\kappa$, and therefore each pole can only be a pole of at most a finite number of terms of the sum \eqref{eq:phi}. This ensures that these points are, in fact, at most simple poles. Let us prove this.

\begin{proposition}\label{pro:finite-fibres}
Given $\sigma_0\in\bR$, there are at most a finite number of multi-indices  $\kappa\in\bN_0^{r-1}$ of length $r-1$ such that $\Re(s_{\kappa})=\sigma_0$.
\end{proposition}

\begin{proof}
The ordering of the roots $\alpha_i$ is irrelevant for our purposes as long as $\alpha_1$ is the root with largest modulus, so we can assume without loss of generality that $\alpha_1>|\alpha_2|\geq\dots\geq |\alpha_r|$. Note that for $\kappa\in\bN_0^{r-1}$,
\[
	(\log\alpha)\:\Re(s_{\kappa})=\kappa_2\log\left|\frac{\alpha_2}{\alpha_1}\right|+\kappa_3\log\left|\frac{\alpha_3}{\alpha_1}\right|+\dots+\kappa_r\log\left|\frac{\alpha_r}{\alpha_1}\right|.
\]

Since $\left|\frac{\alpha_i}{\alpha_1}\right|<1$ for each $2\leq i\leq r$ and $\left|\frac{\alpha_i}{\alpha_1}\right|\geq\left|\frac{\alpha_j}{\alpha_1}\right|$ for $i\leq j$, we have
\[
	|\kappa|\log\left|\frac{\alpha_r}{\alpha_1}\right|\leq \log\alpha\:\Re(s_{\kappa})\leq |\kappa|\log\left|\frac{\alpha_2}{\alpha_1}\right|.
\]
From this, it is clear that
\[
	\Re(s_{\kappa})=\sigma_0\Longrightarrow |\kappa|\leq\frac{\log\alpha}{\log\left|\frac{\alpha_2}{\alpha_1}\right|}\sigma_0,
\]
and this inequality can only be true for a finite number of $\kappa$ (for none, obviously, if $\sigma_0>0$).
\end{proof}

Note that Proposition~\ref{pro:finite-fibres} implies more than the fact that $\iota$ has finite fibres: it also implies that each vertical line can only contain a finite number of ``lines of poles'', and therefore the set $\{s_{\kappa,n}\}$, when varying both $\kappa$ and $n$, is discrete.

%%%%%%%%%%%%%%%%%%%%%%%%%%%%%%%%%%%%%%%%%%%%%%%%%%%%%%%%%%%%%%%%%%%%%%%%%%%%%%%%%%%
\section{Module of relations and the main equivalence}
\label{sec:equivalence}
%%%%%%%%%%%%%%%%%%%%%%%%%%%%%%%%%%%%%%%%%%%%%%%%%%%%%%%%%%%%%%%%%%%%%%%%%%%%%%%%%%%
%%%%%%%%%%%%%%%%%%%%%%%%%%%%%%%%%%%%%%%

We now turn to the matter of determining when the map $\iota$ is injective and, if it is not, study how the points reached by more than one $\kappa\in\bN_0^{r-1}$ behave. It turns out that injectivity is equivalent to a property of diophantine approximation of the number $\alpha$. Following \cite{Dixon} and \cite{Zheng1}, we define:

\begin{definition}\label{d:module-of-relations}
Let $\alpha$ be any algebraic integer of degree $r$, with conjugates $\alpha_1=\alpha,\alpha_2,\dots,\alpha_r$. The \textbf{module of relations} of $\alpha$ is the set
\[
	\cR_{\alpha}=\left\{(m_1,\dots, m_r)\in\bZ^{r}\:\vert\: \alpha_1^{m_1}\dots\alpha_r^{m_r}=1\right\}.
\]
\end{definition}

It is clear that $\cR_{\alpha}$ is a sub-$\bZ$-module of $\bZ^r$. Let us denote by $e=(1,\dots, 1)\in\bZ^r$ the vector all of whose entries are $1$. We can characterize some of the elements of $\cR_{\alpha}$ according to the absolute norm $\norm(\alpha)$ of $\alpha$ by the following statements, with trivial proof:

\begin{lemma}\label{l:trivial-relations}
Let $\alpha$ be any algebraic integer.
\begin{enumerate}[label=(\alph*)]
\item $\norm(\alpha)\neq\pm 1$ if and only if $\cR_{\alpha}\cap\langle e\rangle=\{0\}$.

\item $\norm(\alpha)=-1$ if and only if $\cR_{\alpha}\cap\langle e\rangle=\langle 2e\rangle$.

\item $\norm(\alpha)=1$ if and only if $\cR_{\alpha}\cap\langle e\rangle=\langle e\rangle$.
\end{enumerate}
\end{lemma}

Note that $\cR_{\alpha}=\{0\}$ states the \emph{multiplicative independence} of the algebraic integers $\alpha_1,\dots,\alpha_r$. The condition we are interested in, however, is not multiplicative independence, but \emph{nontrivial} multiplicative independence.

\begin{definition}\label{d:trivial-module}
The module of relations $\cR_{\alpha}$ of $\alpha$ is \emph{trivial} if $\cR_{\alpha}\subseteq\langle e\rangle$, that is, if the only multiplicative relations between the conjugates of $\alpha$ are the ones in which all exponents are equal.
\end{definition}

Note that, if $\alpha$ is a unit, $\cR_{\alpha}$ is trivial if and only if $\rk(\cR_{\alpha})= 1$.

Finally, we denote by $\cH_0$ the hyperplane of $\bZ^r$ ``orthogonal'' to $\langle e\rangle$, that is,
\[
	\cH_0=\left\{(m_1,\dots,m_r\in\bZ^r\:\vert\: m_1+\dots+m_r=0\right\}.
\]

With this, we can now state and prove our main equivalence:

\begin{theorem}\label{t:equivalence}
Let $\alpha$ be a Perron number of degree $r$ with conjugates $\alpha_1=\alpha,\alpha_2,\dots,\alpha_r$. The map $\iota$ defined in \eqref{eq:indexation} is injective if and only if $\cR_{\alpha}\cap\cH_0=\{0\}$.
\end{theorem}

\begin{proof}
For a multi-index $\kappa\in\bN_0^{r-1}$, denote by $\overline{\kappa}\in\bZ^r$ the \emph{extended vector}
\[
	\overline{\kappa}=(-|\kappa|,\kappa_2,\dots,\kappa_r)\in\cH_0.
\]
In addition, given $m=(m_1,\dots,m_r)\in\bZ^{r}$, we denote by $[m]$ the \emph{reduced vector}
\[
	[m]=(m_2,\dots,m_r)\in bZ^{r-1},
\]
so that for a multi-index $\kappa$, $\left[\overline{\kappa}\right]=\kappa$. Note that, for a vector $m\in\bZ^r$, $\overline{[m]}=m\Leftrightarrow m\in\cH_0$.

The theorem follows from the following two statements:

\begin{enumerate}[wide, labelwidth=!, labelindent=0pt, label=(\arabic*)]

\item If $\kappa,\eta\in\bN_0$ and $s_{\kappa}=s_{\eta}$, then $\overline{\kappa-\eta}\in\cR_{\alpha}\cap\cH_0$.

\item Given $m\in\cR_{\alpha}\cap\cH_0$, if $\kappa\in\bN_0$ is such that $\eta=\kappa+[m]\in\bN_0$, then $s_{\kappa}=s_{\eta}$.

\end{enumerate}

Assume that $\kappa,\eta\in\bN_0^{r-1}$ are such that $s_{\kappa}=s_{\kappa}$. Comparing the real parts, we have
\[
	\log\left|(\alpha_1,\dots,\alpha_r)^{\overline{\kappa}}\right|=\log\left|(\alpha_1,\dots,\alpha_r)^{\overline{\eta}}\right|,
\]
which is equivalent to
\[
	\left|(\alpha_1,\dots,\alpha_r)^{\overline{\kappa-\eta}}\right|=1.
\]
Looking now at the imaginary parts, we get that
\[
	\arg\left((\alpha_1,\dots,\alpha_r)^{\overline{\kappa-\eta}}\right)=0\Longrightarrow (\alpha_1,\dots,\alpha_r)^{\overline{\kappa-\eta}}\in\bR_+.
\]

Therefore, $(\alpha_1,\dots,\alpha_r)^{\overline{\kappa-\eta}}$ is a nonnegative real number with absolute value $1$, which means it is $1$ and thus $\overline{\kappa-\eta}\in\cR_{\alpha}\cap\cH_0$. This proves (1), and for (2) the argument can be easily reversed.

From (1) and the observation that $\kappa\neq\eta\Rightarrow\overline{\kappa-\eta}\neq 0$ it follows that $\cR_{\alpha} \cap \cH_{0} \neq \{0\}$ if $\iota$ is not injective.

For the converse, if $m \in \cR_{\alpha} \cap \cH_{0}$ with $m\neq 0$, then $\kappa+[m]\neq\kappa$ for any $\kappa$. Note that there are always $\kappa\in\bN_0^{r-1}$ such that $\kappa+[m]\in\bN_0^{r-1}$ (we can just take $\kappa=(M,\dots, M)$, where $M>\max\{|m_2|,\dots,|m_r|\}$). By (2) this proves that if $\cR_{\alpha} \cap \cH_{0} \neq \{0\}$ then $\iota$ is not injective.
\end{proof}

As it turns out, the condition $\cR_{\alpha}\cap\cH_{0}=\{0\}$ is equivalent, for algebraic units, to the triviality of $\cR_{\alpha}$, by means of the following result:

\begin{lemma}\label{l:rank-formula}
Let $M_1,M_2$ be two sub-$\bZ$-modules of $\bZ^r$. We have that
\begin{equation}\label{eq:rank-formula}
	\rk(M_1)+\rk(M_2)=\rk(M_1+M_2)+\rk(M_1\cap M_2).
\end{equation}
\end{lemma}

\begin{proof}
We have the exact sequence of $\bZ$-modules
\[
	\xymatrix{0\ar[r] & M_1\cap M_2\ar[r] & M_1\times M_2\ar[r] & M_1+M_2\ar[r] & 0}.
\]
Since $\bQ$ is a flat $\bZ$-module and tensoring by $\bQ$ transforms this sequence into a sequence of $\bQ$-vector spaces, converting rank into dimension, we obtain~\eqref{eq:rank-formula}.
\end{proof}

\begin{corollary}\label{c:triviality-intersection}
If $\alpha$ is any algebraic integer of degree $r$ and $\cR_{\alpha}$ is trivial, then $\cR_{\alpha}\cap\cH_0=\{0\}$. Furthermore, this is an equivalence when $\alpha$ is a unit.
\end{corollary}

\begin{proof}
Since $\langle e\rangle\cap\cH_0=\{0\}$, if $\cR_{\alpha}$ is trivial (that is, $\cR_{\alpha}\subseteq\langle e\rangle$), $\cR_{\alpha}\cap\cH_0=\{0\}$.

If $\alpha$ is a unit we have, as we said before, $\cR_{\alpha}$ is trivial if and only if $\rk(\cR_{\alpha})=1$. If $\cR_{\alpha}$ is not trivial, it has rank at least $2$ and has, by the rank formula~\eqref{eq:rank-formula}, nonzero intersection with any hyperplane. In particular, $\cR_{\alpha}\cap\cH_0\neq\{0\}$.
\end{proof}

From this, we can state this following consequence of Theorem~\ref{t:equivalence} and Corollary~\ref{c:triviality-intersection}:

\begin{theorem}\label{t:equivalence-triviality}
Let $\alpha$ be a Perron number. If $\cR_{\alpha}$ is trivial, the pole enumeration map $\iota$ \eqref{eq:indexation} for the zeta function associated to a linear integer recurrence sequence satisfying the minimal polynomial of $\alpha$ is injective. Furthermore, if $\alpha$ is a unit, this is an equivalence.
\end{theorem}

%%%%%%%%%%%%%%%%%%%%%%%%%%%%%%%%%%%%%%%%%%%%%%%%%%%%%%%%%%%%%%%%%%%%%%%%%%%%%%%%%%%
\section{The current state of the question of triviality}
\label{sec:triviality-so-far}
%%%%%%%%%%%%%%%%%%%%%%%%%%%%%%%%%%%%%%%%%%%%%%%%%%%%%%%%%%%%%%%%%%%%%%%%%%%%%%%%%%%
%%%%%%%%%%%%%%%%%%%%%%%%%%%%%%%%%%%%%%%

In the mathematical literature we find many results dealing with conditions that ensure the triviality of the module of relations of an algebraic integer. Most of these give sufficient conditions for triviality, and some are concerned in computing the module of relations, like \cite{Zheng1}. Let us give a brief overview of  sufficiency conditions. Many are stated in terms of the Galois group of the minimal polynomial of $\alpha$. The first one, however, does not: in \cite{Mignotte}, Mignotte states the following.

\begin{theorem}[Mignotte]\label{t:Mignotte}
Let $\theta$ be a Pisot number with conjugates $\theta_1=\theta,\theta_2,\dots,\theta_d$. Then, the relation
\[
	\theta_1^{n_1}\theta_2^{n_2}\dots\theta_d^{n_d}=1
\]
implies $n_1=n_2=\dots=n_d$.
\end{theorem}

Recall that a Pisot number is a real algebraic integer greater than one such that all of its conjugates lie in the open unit disk. In our terms, Mignotte's theorem~\ref{t:Mignotte} states that the module of relations of any Pisot number is trivial.

The second result we cite, which uses Galois theory in its proof but not in its statement, is due to Drmota and Ska\l{}ba. In \cite[Theorem~1]{Drmota1}, they write:

\begin{theorem}[Drmota-Ska\l{}ba]\label{t:Drmota-Skalba}
Let $p$ be an odd prime number, $K$ a field, $f(x)\in K[x]$ an irreducible polynomial (different from $x^p+a_0$) of degree $p$ over $K$ without multiple roots, and $x_1,\dots,x_p\in\overline{K}$ its roots. Furthermore let the degree of the extension $[K(\zeta_m):K]$ be different from $p$ for all positive integers $m$, where $\zeta_m$ denotes a primitive $m$-th root of unity. Then we have for arbitrary integers $k_1,\dots,k_p$
\[
	\prod_{i=1}^p x_i^{k_1}\in K\Longleftrightarrow k_1=k_2=\dots=k_p.
\]
\end{theorem}

When $f(x)\in\bZ[x]$ is the minimal polynomial of an algebraic integer, the conclusion of Theorem~\ref{t:Drmota-Skalba} is once again the triviality of the module of relations. The hypotheses are easily seen to be satisfied when $K=\bQ$ and $f(x)$ is the minimal polynomial of a Perron number of prime degree.

For the general, non-prime, non-Pisot case, we have the following result due to Smyth in \cite[Lemma~1]{Smyth}, which states that numbers with nontrivial module of relations are indeed quite rare:

\begin{theorem}[Smyth]\label{t:Smyth}
Let $\beta$ be an algebraic number, no power of which is rational, with conjugates $\beta=\beta_1,\dots,\beta_n$, and Galois group $\gal(\bQ(\beta)/\bQ)$ the full symmetric group $\cS_n$ on $n$ symbols. Then
\[
	\beta_1^{v_1}\beta_2^{v_2}\dots\beta_n^{v_n}\neq\text{ root of unity}
\]
and
\[
	\sum_{i=1}^{n}\beta_i v_i\neq 0
\]
for any integers $v_1,\dots, v_n$ which are not all equal.
\end{theorem}

This can be generalized to include the alternating groups $A_n$ for $n\geq 4$ by the result \cite[Theorem~1]{Dixon} of Dixon:

\begin{theorem}[Dixon]\label{t:Dixon}
Suppose that $f(X)\in K[X]$ is monic irreducible and the Galois group $G$ of $f(X)$ acts primitively on the set of roots of $f(X)$.

If $G$ acts $2$-transitively on the set of roots of $f(X)$ and there is any nontrivial multiplicative relation between the roots of $f(X)$ then, for some prime $p$, $f(X)$ has degree $p$ and $G$ is isomorphic to the affine group of order $p(p-1)$.
\end{theorem}

An improvement of Dixon's condition of $2$-transitivity was made recently in \cite[Theorem~3.1]{Zheng1}:

\begin{theorem}[Zheng]\label{t:Zheng}
Let $f(x)\in\bQ[x]$ be a univariate polynomial without multiple roots, such that one of its roots is not the root of a rational number. If the Galois group $G_f$, regarded as a permutation group operating on the set of all the complex roots of $f$, is $2$-homogeneous, then the module of relations of $f$ is trivial.
\end{theorem}

Summarising all of these results in our case, we get the following:

\begin{theorem}\label{t:sufficient-conditions}
Let $\alpha$ be a Perron number of degree $r$ and Galois group $G$. If any of the following conditions is satisfied, then the module of relations $\cR_{\alpha}$ is trivial, and therefore the map $\iota$ defined in \eqref{eq:indexation} is injective:
\begin{enumerate}[label=(\arabic*)]
\item $\alpha$ is a Pisot number.

\item $r$ is a prime number.

\item $G$ is the symmetric group $\cS_r$ or the alternating group $A_r$.

\item $G$ acts $2$-homogeneously on the set of conjugates $\alpha_1,\dots,\alpha_r$.
\end{enumerate}
\end{theorem}

\begin{proof}
If $\alpha$ is a Pisot number, the result follows directly from \ref{t:Mignotte}.

If $r$ is an odd prime, it follows from \ref{t:Drmota-Skalba}. If $r=2$, we prove it directly: let $\alpha_1$ be a Perron number of degree $2$, with conjugate $\alpha_2$, and $m_1,m_2\in\bZ$ such that $\alpha_1^{m_1}\alpha_2^{m_2}=1$. By squaring and inverting this relation if necessary, we can assume that $m_1$ and $m_2$ are even and $m_1\geq m_2$. Let $p_0=\alpha_1\alpha_2\in\bZ$. If $m_1=m_2$, we have $p_0=\pm 1$, which implies that $\alpha_1$ is a Pisot number and the result follows from~\ref{t:Mignotte}. If $m_1>m_2$, we have $\alpha_1^{m_1-m_2}=p_0^{-m_2}\in\bQ$, with $m_1-m_2>0$. This means that $\alpha_1$ is a root of the polynomial $x^{m_1-m_2}-p_0^{-m_2}\in\bQ[x]$, but then $\alpha_2$ also would be a root, and this implies that $|\alpha_1|=|\alpha_2|$, contradicting the fact that $\alpha_1$ is a Perron number. Therefore, $\cR_{\alpha}$ is also trivial.

If $G=\cS_r$, the result follows from \ref{t:Smyth}. If it is the alternating group $A_r$, it is $(r-2)$-transitive (see for example \cite[[Theorem~9.7]{Wielandt}), so if $r\geq 4$ we conclude by \ref{t:Dixon}. If $G$ is $A_2$ or $A_3$, $r$ is a prime, the result follows from our previous point.

Finally, if $G$ acts $2$-homogeneously on $\alpha_1,\dots,\alpha_r$, the result follows directly from \ref{t:Zheng}, since a Perron number clearly verifies the other hypotheses.
\end{proof}

%%%%%%%%%%%%%%%%%%%%%%%%%%%%%%%%%%%%%%%%%%%%%%%%%%%%%%%%%%%%%%%%%%%%%%%%%%%%%%%%%%%
\section{Some results on nontriviality}
\label{sec:nontriviality}
%%%%%%%%%%%%%%%%%%%%%%%%%%%%%%%%%%%%%%%%%%%%%%%%%%%%%%%%%%%%%%%%%%%%%%%%%%%%%%%%%%%
%%%%%%%%%%%%%%%%%%%%%%%%%%%%%%%%%%%%%%%

We emphasize that, in spite of Theorem~\ref{t:sufficient-conditions} the module of relations of a Perron number is not always trivial. One example is given by Drmota and Ska\l{}ba in \cite{Drmota1} and \cite{Drmota2}, which they attribute to A. Schinzel: the polynomial
\[
	f(x)=x^6-2x^4-6x^3-2x^2+1
\]
is the minimal polynomial of the Perron number $\alpha=2.24...$ (the image of $\iota$ in this case was shown in Figure~\ref{fig:schinzel}), with conjugates
\[
	\begin{array}{lll}
	\alpha_1=\alpha, & \alpha_3=-0.92...+i1.17..., & \alpha_5=-0.41...+i0.52..., \\
	\alpha_2=0.44..., & \alpha_4=\overline{\alpha_3}, & \alpha_6=\overline{\alpha_5}.\end{array}
\]

The conjugates $\alpha_3$ and $\alpha_4$ line outside the unit disk, so $\alpha$ is not a Pisot number. The module of relations $\cR_{\alpha}$ is nontrivial; in fact, it has rank $4$, generated, for example, by the independent relations
\[
	\alpha_1\alpha_2=1, \quad \alpha_3\alpha_6=1, \quad \alpha_4\alpha_5=1, \quad \alpha_1\alpha_5\alpha_6=1.
\]

Note that the first three of these relations are a result of the fact that $f(x)$ is a reciprocal polynomial, so the roots appear in pairs of inverses, and the fourth comes about because $f(x)$ factors in $\bQ(\sqrt{2})$ as
\[
	f(x)=f_1(x)f_2(x)=(x^3+\sqrt{2}x^2+\sqrt{2}x-1)(x^3-\sqrt{2}x^2-\sqrt{2}x-1),
\]
where $f_1(x)$ has roots $\alpha_2,\alpha_3,\alpha_4$ and $f_2(x)$ has roots $\alpha_1,\alpha_5,\alpha_6$. The Galois group of $f(x)$ is isomorphic to the dihedral group $D_6$. See \cite{Viana} for a deeper study of the Galois groups of reciprocal polynomials.

In general, the module of relations of a number with reciprocal minimal polynomial of degree $r$ (where $r$ must be even, on account of the polynomial being reciprocal) has rank at least $\frac{r}{2}$, arising from the $\frac{r}{2}$ pairings of inverse roots. This would be the case, for example, for all \emph{Salem numbers}. A Salem number is a Perron number such that all of its conjugates lie in the closed unit disk with at least one of them lying on the unit circle. This implies (see for example \cite[Chapter~5]{Pisot-Salem}) that the minimal polynomial of a Salem number of degree $r$ is reciprocal, has exactly two real roots (said Salem number and its inverse, lying inside the open unit disk) and $r-2$ complex roots lying on the unit circle.

Since Salem numbers have reciprocal polynomials, the smallest degree which they can have is $r=4$. The smallest Salem number of degree $4$ is $\beta=1.722...,$, defined as the largest real root of the irreducible polynomial
\[
	g(x)=x^4-x^3-x^2+1,
\]
with Galois group $D_4$. If we denote the conjugates of $\beta$ by
\[
	\begin{array}{ll}
	\beta_1=\beta, & \beta_3=-0.65...+i0.75..., \\
	\beta_2=0.58..., & \beta_4=\overline{\beta_3},\end{array}
\]
the module of relations $\cR_{\beta}$ is generated by the independent relations
\[
	\beta_1\beta_2=1, \quad \beta_3\beta_4=1,
\]
that is, it has rank $2$. Therefore, $\cR_{\beta}\cap\cH_0$ has rank $1$: it is generated by the relation
\[
	\beta_1\beta_2\beta_3^{-1}\beta_4^{-1}=1.
\]

In fact, this is the case for any quartic Salem number:

\begin{proposition}\label{p:Salem-degree-4}
Let $\alpha$ be Salem number of degree $4$. Its module of relations $\cR_{\alpha}$ has rank $2$.

Furthermore, if the conjugates $\alpha_1=\alpha,\alpha_2,\alpha_3,\alpha_4$ are such that $\alpha_2=\alpha_1^{-1}$ and $\alpha_4=\alpha_3^{-1}$, then
\[
	\cR_{\alpha}=\langle (1,1,0,0),(0,0,1,1)\rangle
\]
and
\[
	\cR_{\alpha}\cap\cH_0=\langle (1,1,-1,-1)\rangle.
\]
\end{proposition}

\begin{proof}
The only thing that requires proof is that $\rk(\cR_{\alpha})=2$. Assume
\[
	\alpha_1^{m_1}\alpha_2^{m_2}\alpha_3^{m_3}\alpha_4^{m_4}=1
\]
for some $m_1,...,m_4\in\bZ$. Since $\alpha_2=\alpha_1^{-1}$ and $\alpha_4=\alpha_3^{-1}$, this can be written as
\[
	\alpha_1^{m_1-m_2}\alpha_3^{m_3-m_4}=1.
\]
Taking modules and using the fact that $|\alpha_3|=1$ and $|\alpha_1|>1$, we arrive at $m_1-m_2=0$. From this, $\alpha_3^{m_3-m_4}=1$. However, $\alpha_3$ is not a root of unity, because its minimal polynomial has roots with modulus different than $1$. Therefore, $m_3-m_4=0$, and
\[
	(m_1,m_2,m_3,m_4)\in\langle (1,1,0,0),(0,0,1,1)\rangle.
	\qedhere
\]
\end{proof}

Let us remark on the effect this has on the map $\iota$ in~\eqref{eq:indexation} for any quartic Salem number $\alpha$, with conjugates numbered as in Proposition~\ref{p:Salem-degree-4}. Given $\kappa,\eta\in\bN_0^{3}$, we have that
\[
	s_{\kappa}=s_{\eta} \Longleftrightarrow \overline{\eta}=\overline{\kappa}+m(1,1,-1,-1), \quad m\in\bZ.
\]
Since
\[
	\eta=\left[\overline{\eta}\right]=\left[\overline{\kappa}+m(1,1,-1,-1)\right]=(\kappa_2+m,\kappa_3-m,\kappa_4-m)\in\bN_0^{3},
\]
it must be the case that
\[
	m\geq -\kappa_2, \quad m\leq\kappa_3, \quad m\leq \kappa_4.
\]
This tells us what the size of each fibre is: given $\kappa\in\bN_0^3$,
\[
	\#\iota^{-1}(\{s_{\kappa}\})=\kappa_2+\min\{\kappa_3,\kappa_4\}+1.
\]

In particular, for the pole $s_{(0,0,0)}=0$, the fibre has only one element, which means that there is only one term in the sum \eqref{eq:phi} that has a pole at $s=0$. This is also true for any Perron number, regardless of its module of relations:

\begin{proposition}\label{p:fibre-at-0}
Let $\alpha$ be a Perron number of degree $r$, and $\iota$ defined as in \eqref{eq:indexation}. Then
\[
	\#\iota^{-1}(\{0\})=1.
\]
\end{proposition}

\begin{proof}
Let $\bm{0}=(0,\dots,0)\in\bN_0^{r-1}$, $\iota(\bm{0})=0$. If it were the case that $\#\iota^{-1}(\{0\})>1$, there would be some $\kappa\in\bN_0^{r-1}$, $\kappa\neq\bm{0}$, such that $\iota(\kappa)=0$.

By the same reasoning as in the proof of Theorem~\ref{t:equivalence}, we conclude that $\overline{\kappa}\in\cR_{\alpha}\cap\cH_0$. However, since $\kappa_i\geq 0$ for $i=2,\dots,r$, $\bm{0}+m\kappa\in\bN_0^{r-1}$ for every $m\in\bN$, and again applying Theorem~\ref{t:equivalence} we obtain $\bm{0}+m\kappa\in\iota^{-1}(\{0\})$ for every $m\in\bN_0$. Since not all $\kappa_i$, $i=2,\dots,r$, are $0$, this would contradict the finiteness of the fibre proved in Proposition~\ref{pro:finite-fibres}. Therefore, $\#\iota^{-1}(\{0\})=1$.
\end{proof}

\end{document}